\documentclass[11pt]{article}
\font\smallit=cmti12
\usepackage{amssymb,amsmath,latexsym,amsthm}

\newtheorem{thm}{Theorem}
\newtheorem{lemma}[thm]{Lemma}

\theoremstyle{definition}
\newtheorem{defn}[thm]{Definition}

\makeatletter

\renewcommand\section{\@startsection {section}{1}{\z@}%
{-30pt \@plus -1ex \@minus -.2ex}%
{2.3ex \@plus.2ex}%
{\normalfont\normalsize\bfseries}}

\renewcommand\subsection{\@startsection{subsection}{2}{\z@}%
{-3.25ex\@plus -1ex \@minus -.2ex}%
{1.5ex \@plus .2ex}%
{\normalfont\normalsize\bfseries}}
\renewcommand{\@seccntformat}[1]{\csname the#1\endcsname. } 

\makeatother

\begin{document}

\begin{center}
\uppercase{\bf A Probabilistic Threshold for Monochromatic Arithmetic Progressions}\\[25pt]
{\bf Aaron Robertson} \\
{\smallit Department of Mathematics,
Colgate University, Hamilton, New York}\\
{\tt arobertson@colgate.edu}
\end{center}

\vskip 20pt
\centerline{\bf Abstract}
\noindent
We show that $\sqrt{k}\cdot r^{k/2}$ is a threshold interval length where, under mild conditions, 
almost every $r$-coloring of an interval of longer length contains a monochromatic $k$-term arithmetic progression, while
almost no $r$-coloring of an interval of shorter length contains a
monochromatic $k$-term arithmetic progression.

\baselineskip=15pt

\section{Introduction}

For $k,r \in \mathbb{Z}^+$, let $w(k;r)$ be the minimum integer such that {\it every} $r$-coloring of $[1,w(k;r)]$ admits
a monochromatic $k$-term arithmetic progression.  The existence of such an integer was
shown by van der Waerden \cite{vdw}, and these integers are referred to as van der Waerden
numbers.  Current knowledge places $w(k;r)$ somewhere between
$\frac{r^{k-1}}{ek}(1+o(1))$  and 
$$2^{2{^{r^{2^{2^{k+9}}}}}},$$
with the upper bound being from one of Gowers' seminal works \cite{Gowers}
(with slightly better lower bounds when $r=2$).  A matching of upper and
lower bounds appears unlikely in the near (or distant?) future.  However, by loosening the
restriction that {\it every} $r$-coloring must have a certain property to {\it almost
every} (in the probabilistic sense), we are able to home in on the rate of growth of
the associated numbers.

Although the majority of probabilistic threshold functions studied have concerned
graphs, work on integer Ramsey structures has been investigated.
Different threshold functions related to arithmetic progressions have been studied, in particular, by R\"{o}dl and Ruci\'nski \cite{RR,RR2} and Schacht \cite{Schacht}.

In this article, we assume that every $r$-coloring of a given interval is equally likely, i.e., for each integer in the interval, the probability that
it is a given color is $\frac{1}{r}$.
We refer to a $k$-term arithmetic progression as a
$k$-ap and will use the notation $\langle a,d \rangle_k$ to represent
$a,a+d,a+2d,\dots,a+(k-1)d$, where we refer to $d$ as the {\it gap}
of the $k$-ap.
We  use the standard notation $[1,n]=\{1,2,\dots,n\}$.

\begin{defn} Let $t(k)$ be a function defined on $\mathbb{Z}^+$ with some property $\mathcal{P}$.
We say that
$t(k)$ is a {\it minimal function} (with respect to $\mathcal{P}$) if for every function $s(k)$
defined on $\mathbb{Z}^+$ with property $\mathcal{P}$ we have
$$
\liminf_{k \rightarrow \infty} \frac{t(k)}{s(k)} \leq 1.
$$
We say that
$t(k)$ is a {\it maximal function} (with respect to $\mathcal{P}$) if for every function $s(k)$
defined on $\mathbb{Z}^+$ with property $\mathcal{P}$ we have
$$
\limsup_{k \rightarrow \infty} \frac{t(k)}{s(k)} \geq 1.
$$
\end{defn}

\begin{defn} Let $k,r \in \mathbb{Z}^+$ with $r$ fixed. Denote by $N^+(k;r)$ a minimal function such that
 the probability that a randomly chosen
$r$-coloring of $[1,N^+(k;r)]$  admits a monochromatic $k$-ap tends to 1 as $k \rightarrow \infty$.  Denote by $N^-(k;r)$  a maximal function
such that 
the probability that a randomly chosen
$r$-coloring of $[1,N^-(k;r)]$
admits a monochromatic $k$-ap tends to 0 as $k \rightarrow \infty$. 
\end{defn}

Brown \cite{Brown} showed that $N^+(k;2) \leq (\log k)2^kg(k)$, while Vijay \cite{Vijay}
made a significant improvement by showing that $N^+(k;2) \leq k^{3/2}2^{k/2}g(k)$,
where, in each bound,
$g(k)$ is any function tending to $\infty$.
Vijay \cite{Vijay}, using the linearity of expectation, also provided a lower bound for $N^-(k;2)$ that
is not much smaller than his given upper bound.  The generalization to
$N^-(k;r)$ is straightforward, but included for completeness.

\begin{thm}\label{v} {\rm (Vijay)} Let $f(k) \rightarrow 0$
arbitrarily slowly.  Then $N^-(k;r) \geq \sqrt{k}\cdot r^{k/2}f(k)$.
\end{thm}

\begin{proof} Let $n=\sqrt{k}\cdot r^{k/2}f(k)$.  We will show that
the probability that a random $r$-coloring of $[1,n]$ admits a
monochromatic $k$-ap tends to $0$ (as $k \rightarrow \infty$) by calculating
the expected number of monochromatic $k$-aps in two ways.
Consider a random $r$-coloring $\chi$ of $[1,n]$.
Let $X_i$ equal $1$ if the $i^{\mathrm{th}}$ $k$-ap is monochromatic
under $\chi$; otherwise, let $X_i$ equal $0$.  Then
$$
X=\sum_{i=1}^{\frac{n^2}{2(k-1)}}  X_i
$$
is the number of monochromatic $k$-aps under $\chi$ (where we are suppressing the asymptotic $(1+o(1))$ in the upper limit).  By the
linearity of expectation, we have
$$
E(X) = \sum_{i=1}^{\frac{n^2}{2(k-1)}} E(X_i) = \sum_{i=1}^{\frac{n^2}{2(k-1)}} P(X_i=1) =  \frac{n^2}{2(k-1)} \cdot \frac{1}{r^{k-1}}
= \frac{rk}{2(k-1)} \cdot f^2(k).
$$

Using the definition of $E(X)$, we have
$$
E(X) = \sum_{i=0}^{\frac{n^2}{2(k-1)}} i \cdot P(X=i)
\geq \sum_{i=0}^{\frac{n^2}{2(k-1)}} P(X=i) - P(X=0)
= 1-P(X=0).
$$
Hence, we have
$
\frac{rk}{2(k-1)} \cdot f^2(k) \geq 1-P(X=0), 
$ so that
$$
P(X=0) \geq 1-\frac{rk}{2(k-1)} \cdot f^2(k) \underset{k \rightarrow \infty}{\longrightarrow} 1,
$$
which completes the proof.
\end{proof}

We note here that if we consider the set of $k$-aps with gaps that are primes larger than $k$, we can
follow Vijay's argument for an upper bound on $N^+(k;2)$ very closely to match his bound: $N^+(k) \leq k^{3/2}\,2^{k/2}
g(k)$ for any function $g(k) \rightarrow \infty$  (we will use
the notation $g(k) \rightarrow \infty$ as
opposed to $g(k) \rightarrow 0$ as found in \cite{Vijay}).  In the next section,
we construct a larger family of $k$-aps (than Vijay's and than those $k$-aps with prime gap larger than $k$) with the aim of lowering this upper bound 
on $N^+(k;2)$ by a
factor of $k$, while also generalizing to an arbitrary number of colors.

\section{A Structured Family of Arithmetic Progressions}

For $k,n\in \mathbb{Z}^+$, let
$
AP_k(n) = \{\langle a,d \rangle_k \subseteq [1,n]: a,d \in \mathbb{Z}^+\},
$
i.e., the set of $k$-aps in $[1,n]$  and let $A_j(n)=\{\langle a,d \rangle_k \in AP_k(n): d=j\}$
so that $AP_k(n)$ is the disjoint union of the
$A_j(n)$:
$$
AP_k(n) = \bigsqcup_{d=1}^{\left\lfloor\frac{n-1}{k-1}\right\rfloor} A_d(n).
$$

We will now sieve out elements from each $A_d(n)$.  
Order the elements of $A_d(n)$ by their initial term.
Remove every $k$-ap in $A_d(n)$ with initial term in
$$
\bigcup_{j=1}^{\left\lceil\frac{n}{3d}-\frac{k-1}{3}\right\rceil} \big[(3j-2)d+1,3jd\big].
$$
Let $\overline{A}_d(n)$ be the set of elements in $A_d(n)$ which
were not removed and define
\begin{equation*}
\overline{AP}_k(n) = \bigsqcup_{d=1}^{\left\lfloor\frac{n-1}{k-1}\right\rfloor} \overline{A}_d(n).
\end{equation*}
Notice that, for each $d$, every two {\it intersecting} $k$-aps in $\overline{A}_d(n)$ have
initial terms that are more than $2d$ apart, by construction.

We next provide a useful lemma which gives bounds on the size of 
$\overline{AP}_k(n)$.

\begin{lemma}\label{lem1} Let $k\geq 5$. Then
\begin{equation}\label{eqn:2}
\frac{|AP_k(n)|}{3} \leq |\overline{AP}_k(n)| \leq \frac{|AP_k(n)|}{2}.
\end{equation}
Hence,
\begin{equation}\label{eqn:3}
\frac{n^2}{6(k-1)}(1+o(1)) \leq \left|\overline{AP}_k(n)\right| \leq 
\frac{n^2}{4(k-1)}(1+o(1)).
\end{equation}
\end{lemma}

\begin{proof}
It is a standard exercise to show that $\left|{AP}_k(n)\right|=\frac{n^2}{2(k-1)}(1+o(1))$, so that the inequalities in (\ref{eqn:3})
follow immediately once we prove the inequalities in (\ref{eqn:2}).
The lower bound in (\ref{eqn:2}) follows easily from the
sieving construction, hence we focus on the upper bound in (\ref{eqn:2}).

For $d \in \left[1, \frac{n}{k+6}\right)$ we have $\overline{A}_d(n) \leq 
\frac37 |A_d(n)|$. To see this, note that since $d<\frac{n}{k+6}$ we
have $7d + (k-1)d \leq n$, so that $A_d(n)$ has at least $7d$ elements.
Of the $k$-aps in $A_d(n)$ with initial term at most $7d$,
our sieve admits exactly $\frac37^{\mathrm{ths}}$ of them to 
$\overline{A}_d(n)$. By continuing the sieving process, starting with
the $k$-ap with initial term $7d+1$ (if it exists), we will remove
at least twice as many $k$-aps from $A_d(n)$ than we add to $\overline{A}_d(n)$;
hence, $\overline{A}_d(n) \leq 
\frac37 |A_d(n)|$ holds.

For $d \in \left[\frac{n}{k+6} ,\frac{n-1}{k-1}\right]$ we use the trivial
bound $\overline{A}_d(n) \leq A_d(n)$. 

Lastly, we note that $|A_d(n)|= n-(k-1)d$.
Hence,
$$
\begin{array}{rl}
\displaystyle\left|\overline{AP}_k(n)\right|
= \sum_{d=1}^{\left\lfloor\frac{n-1}{k-1}\right\rfloor} |\overline{A}_d(n)| &\displaystyle = \sum_{d=1}^{\left\lfloor\frac{n}{k+6}\right\rfloor - 1} |\overline{A}_d(n)| + \sum_{d=\left\lfloor\frac{n}{k+6}\right\rfloor}^{\left\lfloor\frac{n-1}{k-1}\right\rfloor} |\overline{A}_d(n)|\\[25pt]
&\displaystyle \leq 
\sum_{d=1}^{\left\lfloor\frac{n}{k+6}\right\rfloor - 1} \frac{3|A_d(n)|}{7} + \sum_{d=\left\lfloor\frac{n}{k+6}\right\rfloor}^{\left\lfloor\frac{n-1}{k-1}\right\rfloor} |A_d(n)|\\[25pt]
&\displaystyle=\frac{3k^3+33k^2+23k+284}{14(k-1)^2(k+6)^2}\cdot n^2(1+o(1))\\[20pt]
&\displaystyle \leq \frac{1}{4(k-1)}n^2 (1+o(1)) \qquad \mbox{\small (for $k \geq 5$)}\\[20pt]
&\displaystyle = \frac{|AP_k(n)|}{2},
\end{array}
$$
which proves the upper bound in (\ref{eqn:3}).
\end{proof}

We next give a result on how pairs of $k$-aps from $\overline{AP}_k(n)$
intersect.

\begin{lemma}\label{lem:4} Let $A=\langle a,b \rangle_k$ and $C=\langle c,d \rangle_k$ belong to $\overline{AP}_k(n)$. 
Then $|A\cap C| \leq k-3$.  Furthermore,
\begin{itemize}
\item[(i)] $|A \cap C| > \left\lceil\frac{k}{2}\right\rceil$ only if $b=d$
\item[(ii)] $\left\lceil\frac{k}{3}\right\rceil\leq|A \cap C|\leq\left\lceil\frac{k}{2}\right\rceil$ only if
$\frac{b}{d}\in \left\{\frac13, \frac12,\frac23,1,\frac32,2,3\right\}$.
\end{itemize}
\end{lemma}

\noindent
{\it Proof.} 
We first argue that in order to have $|A \cap C| \geq \left\lceil\frac{k}{3}\right\rceil$
we must have $\frac{b}{d}\in \left\{\frac13, \frac12,\frac23,1,\frac32,2,3\right\}$.
Consider $b \leq d$. 
We have
$a+ib=c+j_1d$ and $a+(i+x)b=c+j_2d$ for some $i\in[0,k-2]$ and $x \in \{1,2,3\}$ (else we cannot have enough intersections) and $j_1<j_2$.
Thus, $xb=(j_2-j_1)d$.  Since $d\geq b$ we must have $j_2-j_1\leq x$. 
This leaves $\frac{b}{d} = \frac{j_2-j_1}{x} \in \left\{\frac13,\frac12,\frac23,1\right\}$.
For $d > b$, we take reciprocals and achieve the stated goal.

Now, in order for $A$ and $C$ to intersect in more  than $\left\lceil\frac{k}{2}\right\rceil$ places, there
must be two consecutive elements of, say, $A$ in the intersection.
Let $a+ib$ and $a+(i+1)b$ be two such elements.  We must have $d \leq b$ in order for
$C$ to intersect both of these.  So, let $a+ib = c+jd$ and $a+(i+1)b = c+\ell d$.
These imply that $b = (\ell-j)d$.  If $\ell-j > 1$ then $d \leq \frac{b}{2}$.  In this
situation, $C$ intersects $A$ in at most $\left\lceil\frac{k}{2}\right\rceil$ places since for
every two consecutive terms of $A$, there exists a term of $C$ between them.
Thus, $\ell-j=1$ and $b=d$ as stated.

To show that $|A \cap C| \leq k-3$, note that we have proved part (i) so we need only
consider $k$-aps with the same gap, i.e., those in the same $\overline{A}_g(n)$ for some gap $g$.  In order for
two such $k$-aps to intersect in more than $k-3$ places, their starting elements must be within
$2g$ of each other.  As noted before, by construction of $\overline{A}_g(n)$, this is not possible.
\hfill $\Box$

\begin{lemma}\label{lem6}
For a given $A=\langle a,b\rangle_k\in \overline{AP}_k(n)$, the number of $\langle c,d\rangle_k\in \overline{AP}_k(n)$ with $c\geq a$ that intersect $A$ in $p$ places is 
\begin{itemize}
\item[\it (i)] $0$ for $p>k-3$;
\item[\it (ii)] $1$ for each $p \in \left[ \left\lceil\frac{k}{2}\right\rceil+1,k-3\right]$;
\item[\it (iii)] at most $7$ for each $p\in \left[ \left\lceil\frac{k}{3}\right\rceil,\left\lceil\frac{k}{2}\right\rceil\right]$.
\end{itemize}
\end{lemma}

\noindent
{\it Proof.} Part (i) is just a restatement of part of Lemma \ref{lem:4}.  For part (ii), by Lemma \ref{lem:4}(i), we
must have $b=d$.  For a given $p$, we have $c=a+(k-p)d$ and the result follows. For part (iii),
since $b$ is fixed, $d$ must be one of the 7 gaps that adhere to Lemma \ref{lem:4}(ii).  In order to
intersect in exactly $p$ places, $c$ is determined.
\hfill $\Box$

\vskip 5pt
With these lemmas under our belt, we are now ready to move onto the main result.

\section{Main Result}

We incorporate Theorem \ref{v} into the main result, which we now state.

\begin{thm} \label{thm:7}
Let $f(k) \rightarrow 0$ and $g(k) \rightarrow \infty$ arbitrarily slowly.  Then,
$$
\sqrt{k}\cdot r^{k/2}f(k) \leq N^-(k;r)< N^+(k;r) \leq \sqrt{k}\cdot r^{k/2}g(k).
$$
\end{thm}

\begin{proof}  We need only to prove the upper bound on $N^+(k;r)$ and do so by using the family
defined in the previous section, along with techniques from \cite{Brown} and \cite{Vijay}.  To this end,
let $n = \sqrt{k}\cdot r^{k/2}g(k)$
and partition $[1,n]$ into intervals of length $s = \left\lceil\frac{n}{g(k)^{4/3}}\right\rceil$,
where the last interval may be shorter (and won't be used).
Let $I_j$, $1 \leq j \leq \lfloor \frac{n}{s}\rfloor$ be these intervals.  

Fix an interval $I_\ell$ and randomly color each element of $I_\ell$ with one
of $r$ colors, where each color is equally likely.
Let $X_i$ be the event that the $i^{\mathrm{th}}$
$k$-ap in $I_\ell$ is monochromatic, where $1 \leq i \leq \frac{s^2}{6(k-1)}$ holds
by the lower bound in (\ref{eqn:3}) (and we suppress lower order terms).
Denote by $p$  the probability that a
 monochromatic $k$-ap exists. Via one of the Bonferroni inequalities, we have
$$
p=\left|P\left(\bigcup_{i=1}^{\frac{s^2}{6(k-1)}} X_i\right)\right| \geq \sum_{i=1}^{\frac{s^2}{6(k-1)}} P(X_i) \hskip 10pt- \hskip -5pt\sum_{1 \leq i<j \leq \frac{s^2}{6(k-1)}} \hskip -10pt P(X_i \cap X_j).
$$
Hence,
$$
p \geq \frac{s^2}{6(k-1)} \cdot \frac{1}{r^{k-1}} \hskip 10pt- \hskip -5pt \sum_{1 \leq i<j \leq \frac{s^2}{6(k-1)}} \hskip -10pt P(X_i \cap X_j).
$$

We now focus on the double summation.  With a slight abuse of notation, we rewrite this as
$$
\sum_{b\in\overline{AP}_k(s)}\sum_{a\in\overline{AP}_k(s)}P(X_a \cap X_b),
$$
where the initial term of $b$ is at most as large as the initial term of $a$.  For
a given $b \in \overline{AP}_k(s)$ with gap $g$, we define:
\vskip 5pt
\hskip 60pt$S_b = \{a\in\overline{A}_g(s): a \cap b \neq \emptyset\};$

\vskip 5pt
\hskip 60pt$T_b =\left\{a \in \overline{A}_h(s): \left\lceil\frac{k}{3}\right\rceil
\leq |a \cap b| \leq \left\lceil\frac{k}{2}\right\rceil \mbox{ and } \frac{g}{h} \in \{\frac13,\frac12,\frac23,\frac32,2,3\}\right\}$ 

\vskip 5pt
\hskip 60pt$Q_b= \left\{a \in \overline{AP}_k(s): a\cap b \neq \emptyset \mbox{ and } a \not \in S_b \sqcup T_b\right\};$
\vskip 5pt
\hskip 60pt$R_b =\{a \in \overline{AP}_k(s): a\cap b = \emptyset\}$.
\vskip 5pt
\noindent
Note that $\overline{AP}_k(s)$ is the disjoint union $S_b \sqcup T_b \sqcup Q_b \sqcup R_b$.
   With these definitions, the double summation becomes
$$
\sum_{b\in\overline{AP}_k(s)}\left(
\sum_{a \in S_b} P(X_a \cap X_b)
+
\sum_{a \in T_b} P(X_a \cap X_b)
+
\sum_{a \in Q_b} P(X_a \cap X_b)
+
\sum_{a \in R_b} P(X_a)P(X_b)
\right).
$$

Appealing to Lemmas \ref{lem1} and \ref{lem:4}, we find the following (for $k$ sufficiently large):
\begin{equation}\label{samegap}
\sum_{a \in S_b} P(X_a \cap X_b) \leq \sum_{i=2}^{k-1} \frac{1}{r^{k+i}} = \frac{1}{r^{k+1}} - \frac{1}{r^{2k-1}}
\leq \frac{1}{r^{k+1}};
\end{equation}
\begin{equation}\label{gcd23}
\hskip -19pt\sum_{a \in T_b} P(X_a \cap X_b) \leq \frac{7\left(\left\lceil\frac{k}{2}\right\rceil - \left\lceil\frac{k}{3}\right\rceil+1\right)}{r^{k+k/2}}\leq \frac{1}{3\cdot r^{k+1}};
\end{equation}
\begin{equation}\label{Q}
\hskip -70pt\sum_{a \in Q_b} P(X_a \cap X_b) \leq \frac{3sk}{r^{k+2k/3}}\leq \frac{1}{3\cdot r^{k+1}};
\end{equation}
\begin{equation}\label{R}
\hskip -40pt\sum_{a \in R_b} P(X_a \cap X_b) \leq \frac{s^2}{4(k-1)r^{2k-2}}\leq \frac{1}{3\cdot r^{k+1}},
\end{equation}
where: \eqref{samegap} holds since there is exactly one such $k$-ap that intersects $b$
in $k-1-i$ places and Lemma \ref{lem6}(i) gives $i \geq 2$;
\eqref{gcd23} follows from Lemma \ref{lem:4}(ii) and Lemma \ref{lem6}(iii);
\eqref{Q} holds from the lower bound in Lemma \ref{lem:4}(ii) and since $3sk$ or fewer $k$-aps intersect a given $k$-ap (this is a standard bound
typically used with the Lovasz Local Lemma; see, e.g., \cite{BLR}); and
\eqref{R} holds by using the upper bound in (\ref{eqn:3}) along with independence since the two $k$-aps do not share any element.

Using these bounds, we have
$$
\sum_{a \in S_b} P(X_a \cap X_b)
+
\sum_{a \in T_b} P(X_a \cap X_b)
+
\sum_{a \in Q_b} P(X_a \cap X_b)
+
\sum_{a \in R_b} P(X_a \cap X_b)
\leq \frac{2}{r^{k+1}} \leq \frac{1}{r^k}.
$$
Hence, using the upper bound in (\ref{eqn:3}), we have
$$
\sum_{b\in\overline{AP}_k(s)}\sum_{a\in\overline{AP}_k(s)}P(X_a \cap X_b) \leq \frac{s^2}{4(k-1)} \cdot \frac{1}{r^{k}}
$$
so that
$$
p \geq \frac{s^2}{6(k-1)} \cdot \frac{1}{r^{k-1}} - \frac{s^2}{4(k-1)} \cdot \frac{1}{r^{k}}
= \frac{(2r-3)s^2}{12(k-1)} \cdot \frac{1}{r^{k}}
$$

This gives us that the probability that $I_\ell$ has no monochromatic
$k$-ap from $\overline{AP}_k(s)$ is at most
$$
1 - \frac{(2r-3)s^2}{12(k-1)r^k}.
$$
Thus, the probability that none of the intervals $I_j$ contains a monochromatic $k$-ap from $\overline{AP}_k(n)$ is,
for $k$ sufficiently large,
at most
$$
\begin{array}{rl}
\left(
1 - \frac{(2r-3)s^2}{12(k-1)r^k}
\right)^{g^{4/3}(k)-1} \lessapprox \exp\left(-\frac{(g^{4/3}(k)-1)(2r-3)s^2}{12(k-1)r^k}\right)
&\leq \exp\left(-\frac{(2r-3)n^2}{g^{5/3}(k)\cdot12(k-1)r^k}\right)\\[15pt]
&=\exp\left(-\frac{k(2r-3)}{12(k-1)}g^{1/3}(k)\right) \underset{k \rightarrow \infty}{\longrightarrow} 0.
\end{array}
$$
It follows that
the probability that one of the intervals $I_j$ 
contains a monochromatic $k$-ap from $\overline{AP}_k(n)$ tends to 1.  Hence,
the probability that a random $r$-coloring of $[1,n]$ admits a monochromatic
$k$-ap also tends to 1, thereby completing the proof.
\end{proof}

Although  Theorem \ref{thm:7} shows that $\sqrt{k}\cdot r^{k/2}$ is a
threshold function, we wonder whether $f(k)$ and $g(k)$ 
can be replaced by constants as results in \cite{RR,RR2,Schacht} have.


\begin{thebibliography}{9} \footnotesize \parskip=1pt

\bibitem{Brown} Tom Brown, A pseudo upper bound for the van der Waerden function,
{\it J. Combin. Theory Ser. A} {\bf 87} (1999), 233-238.

\bibitem{BLR}  Tom Brown, Bruce Landman, and Aaron Robertson, Bounds on some van der Waerden numbers,
{\it J. Combin. Theory Ser. A} {\bf 115}, 1304-1309.

\bibitem{Gowers} W. T. Gowers, A new proof of Szemer\'{e}di's theorem, {\it Geom. Funct.
Anal.} {\bf 11} (2001), no. 3, 465-588.

\bibitem{RR}  Vojt\v{e}ch R\"odl and Andrzej Ruci\'nski, Threshold functions for
Ramsey properties, {\it J. Amer. Math. Soc.} {\bf 8} (1995), 917-942.

\bibitem{RR2}  Vojt\v{e}ch R\"odl and Andrzej Ruci\'nski,
Rado partition theorem for random subsets of integers,
{\it Proc. London Math. Soc. (3)} {\bf 74} (1997), 481-502.

\bibitem{Schacht} Mathias Schacht, Extremal results for random
discrete structures, preprint.

\bibitem{vdw} B. L. Van der Waerden, Beweis einer baudetschen Vermutung,
{\it Nieuw Archief voor Wiskunde} {\bf 15} (1927), 212-216.

\bibitem{Vijay} Sujith Vijay, Monochromatic progressions in random colorings, {\it J. Combin. Theory Ser. A}
{\bf 119} (2012), 1078-1080.





\end{thebibliography}
 \end{document}